# On Rings of MAL'CEV–NEUMANN Series

Mohammad. H. Fahmy, Refaat. M. Salem and

Shaimaa. Sh. Shehata

Department of Mathematics, Faculty of Science,

Al-Azhar University,

Nasr City, 11884, Cairo, Egypt.

Dedicated to the memory of Professor Mohammad H. Fahmy

Abstract. In this paper, we investigate the conditions for the Mal'cev–Neumann series ring $\Lambda = R((G; \sigma; \tau))$ to be left fusible and an $SA$-ring. Also, we show that: if $G$ is a quasitotally ordered group and $U$ a $\Sigma$-compatible semiprime ideal of $R$, then $R((G; \sigma; \tau))$ is a $\Sigma_{U((G;\sigma;\tau))}$-zip ring if and only if $R$ is a $\Sigma_U$-zip ring.

**MSC 2020:** 16D25, 16U80, 16U99.

**Keywords:** Σ-zip ring, fusible rings, Mal'cev-Neumann series.

## 1. Introduction

Throughout this paper $R$ denotes an associative ring with identity, $r_R(X)$ the right annihilator of $X$ in $R$ for any subset $X$ of a ring $R$, $nil(R)$ the set of all nilpotent elements of $R$ and for any two nonempty subsets $U$ and $V$ of $R$, let $(U:V) = \{x \in R \mid Vx \subseteq U\}$. It is easy to see that if $U$ and $V$ are two right ideals of $R$, then $(U:V)$ is an ideal of $R$ and such ideal is usually called the quotient of $U$ by $V$. Section 2, is devoted to recall some background concerning the structure of the ring $\Lambda = R((G; \sigma; \tau))$ of Mal'cev–Neumann series. In section 3, we focus on a property of nonzero zero divisor elements related to the fact that the sum of two zero-divisors need not be a zero-divisor. So, the set of left zero-divisors in a ring R is not a left ideal. Therefore, there exists a left zero-divisor which can be expressed as the sum of a left zero-divisor and a non-left zero-divisor in R. This leads Ghashghaei & McGovern [7] to introduce a class of rings in which every element can be written as the sum of a left zero-divisor and left regular element. They called this class of rings left fusible. This leads us to extend the left fusible property of $R[x, \sigma]$ and $R[\![x, \sigma]\!]$ [7, Proposition 2.9] to the ring $\Lambda = R((G; \sigma; \tau))$ of Mal'cev–Neumann series in Proposition 3.2. In section 4, we discuss a class of rings introduced by Birkenmeier et al [2] called right $SA$-ring, it is exactly the class of rings for which the lattice of right annihilator ideals is a sub lattice of the lattice of ideals. This class includes all quasi-Baer (hence all Baer) rings and all right $IN$-rings (hence all right self injective rings). They showed that this class is closed under direct products, full and upper triangular matrix rings, certain polynomial rings, and two-sided rings of quotients. This drives us to extend the SA property of $R[x]$ [2, Theorem 3.2] to the ring $\Lambda = R((G; \sigma; \tau))$ of Mal'cev–Neumann series in Theorem 4.5. In section 5, we discuss the class of right zip rings introduced by Faith [6] and its generalizations. A ring $R$ is called right zip provided that if the right annihilator $r_R(X)$ of a subset $X$ of $R$ is zero, then there exists a finite subset $Y \subseteq X$ such that $r_R(Y) = 0$. $R$ is zip if it is both right and left zip. The concept of zip rings was initiated by Zelmanowitz [18] and appeared in enormous papers [see,1, 4, 5, 6, 8 and 15] and references therein. Then Ouyang [13] generalized this concept through the introduction of the notion of right weak zip rings (i.e., rings provided that if the right weak annihilator of a subset $X$ of $R$, $N_R(X) \subseteq nil(R)$, then there exists a finite subset $X_0 \subseteq X$ such that $N_R(X_0) \subseteq nil(R)$, where $N_R(X) = \{a \in R \mid xa \in nil(R) \text{ for each } x \in X\}$. Ouyang studied the transfer of the right (left) weak zip property between the base ring $R$ and some of its extensions such as the ring of upper triangular matrices $T_n(R)$ and Ore extension $R[x, \sigma, \delta]$, where $\sigma$ is an endomorphism and $\delta$ is a $\sigma$-derivation. Also, in [14] Ouyang et al continued the study of zip rings and introduced the notion of a $\Sigma$-zip ring and investigated some of its properties. For a proper ideal $U$ of $R$, $R$ is called $\Sigma_U$-zip provided that for any subset $X$ of $R$ with $X \nsubseteq U$, if $(U:X) = U$, then there exists a finite subset $Y \subseteq X$ such that $(U:Y) = U$. Clearly, if $U = 0$, then for any subset $X$ of $R$, we have $(U:X) = r_R(X)$, and so $R$ is $\Sigma_0$-zip if and only if $R$ is right zip. If $R$ is an $NI$ ring (i.e., a ring in which $nil(R)$ forms an ideal) and $U = nil(R)$. Then for any subset $X$ of $R$, we have $(nil(R):X) = N_R(X)$, and so $R$ is $\Sigma_{nil(R)}$-zip if and only if $R$ is weak zip. So, both right zip rings and weak zip rings are special cases of $\Sigma$-zip rings. This caused us to pay attention to prove that, $R$ is $\Sigma_U$-zip ring if and only if $\Lambda = R((G;\sigma;\tau))$ is a $\Sigma_{U((G;\sigma;\tau))}$-zip ring.

## 2. Rings of Mal'cev–Neumann Series

The Mal'cev-Neumann construction appeared for the first time in the latter part of 1940's (the Laurent series ring, a particular case of Mal'cev-Neumann rings, was used before by Hilbert). Using them, Mal'cev and Neumann independently showed (in 1948 and 1949 resp.) that the group ring of an ordered group over a division ring can be embedded in a division ring. Since then, the construction has appeared in many papers, mainly in the study of various properties of division rings and related topics. For in-stance, Makar-Limanov in [10] used a skew Laurent series division ring to prove that the skew field of fractions of the first Weyl-algebra contains a free noncommutative subalgebra. Other results on Mal'cev-Neumann rings can be found in Lorenz [9], Musson and Stafford [11], Sonin [17] and Zhao and Liu. [19]. A pair $(G, \leq)$, where $G$ is a group and $\leq$ an order relation, is called quasitotally if the order $\leq$ can be refined to a total order $\preccurlyeq$ on $G$. Let $(G, \leq)$ be a quasitotally ordered group, and $\sigma$ a map from $G$ into the group of automorphisms of $R(Aut(R))$, which assigns to each $x \in G, \sigma_x \in Aut(R)$. Suppose also that we are given a map $\tau$ from $G \times G$ to $U(R)$, the group of invertible elements of $R$. Let $\Lambda = R((G;\sigma;\tau))$ denotes the set of all formal sums $f = \sum_{x \in G} r_x \bar{x}$ with $r_x \in R$ such that $supp(f) = \{ x \in G \mid r_x \neq 0\}$ (the support of $f$) is a artinian and narrow subset of $G$, with component wise addition and multiplication defined by: $(\sum_{x \in G} a_x \bar{x})(\sum_{y \in G} b_y \bar{y}) = \sum_{z \in G} \left( \sum_{\{x,y \mid xy = z\}} a_x \sigma_x(b_y) \tau(x,y) \right) \bar{z}$. In order to insure associativity of $\Lambda$, it is necessary to impose two additional conditions on $\sigma$ and $\tau$, namely that for all $x, y, z \in G$,

(i)   $\tau(xy; z)\sigma_x(\tau(x; y)) = \tau(x; yz)\tau(y; z)$;

(ii)  $\sigma_y \sigma_z = \sigma_{yz} \eta(y; z)$;

where $\eta(y; z)$ denotes the automorphism of $R$ induced by the unit $\tau(y; z)$. It is now routine to check that $\Lambda = R((G; \sigma; \tau))$ is a ring which is called the Mal'cev-Neumann series ring, that has as an $R$-basis, the set $\bar{G}$ (a copy of $G$) and contains $R$ as a subring where the embedding of $R$ into $\Lambda$ is given by $r \to r\bar{1}$. It is easy to see that the identity element of $\Lambda$ is of the form $1 = u\bar{1}$ for some $u \in U(R)$, unlike group rings, crossed product also rings of Mal'cev–Neumann series do not have a natural basis. Indeed if $d: G \to U(R)$ assigns to each element $g \in G$ a unit $d_g$, then $\widetilde{G} = \{\tilde{g} = g d_g | g \in G\}$ yields an alternative $R$-basis for $\Lambda$ which still exhibit the basic Mal'cev–Neumann structure and this is called a diagonal change of basis. Thus, via diagonal change of basis we can still assume that $1 = \bar{1}$.

In [14] Ouyang called an ideal $U$ semiprime if for any $a \in R, a^n \in U$ implies $a \in U$. We denote $U((G; \sigma; \tau))$ the subset of $\Lambda$ consisting of those elements whose coefficients lie in $U$, that is, $U((G; \sigma; \tau)) = \{f = \sum_{x \in G} a_x x \in R((G; \sigma; \tau)) | a_x \in U, x \in suppf\}$ For each $f \in \Lambda$, let $C(f)$ be the content of $f$, i.e. $C(f) = \{a_x \mid x \in supp(f)\}$.

For $f = (\sum_{x \in G} a_x x)$ and $g = (\sum_{y \in G} b_y y) \in \Lambda$ we define $X_w(f, g) = \{(x_i, y_j) \in G \times G | x_i y_j = w$ where $x_i \in supp(f)$ and $y_j \in supp(g)\}$. It is well known that $X_w(f, g)$ is a finite subset. Let $R$ be a ring and $G$ a quasitotally ordered group, $R$ is called a $G$-Armendariz ring if whenever $f = \sum_{x \in G} a_x x$ and $g = \sum_{y \in G} b_y y \in \Lambda = R((G; \sigma; \tau))$ such that $fg = 0$ implies $a_x b_y = 0$ for each $x \in supp(f)$ and $y \in supp(g)$.

## 3. Fusible rings of Mal'cev–Neumann Series.

It is well known that an element $a \in R$ is a left zero-divisor if there is $0 \neq r \in R$ with $ar = 0$ and an element which is not a left zero-divisor is called a non-left zero-divisor. An element $a \in R$ is regular if it is neither a left zero-divisor nor a right zero-divisor. Let $Z_\ell(R)$ (respectively, $Z_\ell^*(R)$) denote the set of left zero-divisors (respectively, non-left zero-divisors) of $R$. Similarly, $Z_r(R)$ (respectively, $Z_r^*(R)$) denote the set of right zero-divisors (respectively, non-right zero-divisors) of $R$.

In this section, we study the left fusible and right nonsingular rings of Mal'cev–Neumann series.

**Definition 3.1.** A ring R is said to be left fusible if every element can be expressed as a sum of a left zero divisor and a non-left zero divisor (left regular).

**Proposition 3.2.** Let $(G, \leq)$ be a quasitotally ordered group, $\sigma: G \to Aut(R)$ and $\tau: G \times G \to U(R)$, the group of units in $R$. If $R$ is a $\sigma$- compatible and left fusible ring, then $R((G; \sigma; \tau))$ is left fusible.

**Proof.** Let $R$ be a left fusible ring and $0 \neq f \in R\big((G; \sigma; \tau)\big)$. Since, the order $\leq$ on $G$ can be refined to a total ordered $\leqslant$ on G, then there exists $0 \neq s_0 = \pi(f) \in G$ a minimal element in $supp\ f$ with respect to $\leqslant$. Since $R$ is a left fusible ring, then there exists $a \in Z_\ell(R)$ and $b \in Z_\ell^*(R)$ such that $f(s_0) = a + b$. Since $a \in Z_\ell(R)$, then there exists $d \in R$ such that $ad = 0$. Now consider $g, h \in R\big((G; \sigma; \tau)\big)$ such that

$$g(s) = \begin{cases} a & \text{if } s = s_0 \\ 0 & \text{otherwise} \end{cases} \text{ and } h(s) = \begin{cases} b & \text{if } s = s_0 \\ f(s) & \text{if } s \neq s_0 \end{cases}, \text{ consequently } f = g + h.$$

Since $R$ is $\sigma$-compatible and $\tau(1, s_0)$ is an invertible element it follows that $0 = ad = g(s_0)d = g(s_0)\sigma_{s_0}(d) = g(s_0)\sigma_{s_0}(d)\tau(1, s_0) = (gd\bar{1})(s_0)$. Hence $g \in Z_\ell(R((G; \sigma; \tau)))$. Now we prove that $h \in Z_\ell^*(R((G; \sigma; \tau)))$. To the contrary suppose that $h \in Z_\ell(R((G;\sigma;\tau)))$, so there exists $k \in R((G;\sigma;\tau))$ such that $hk = 0$. By hypothesis the order $\leq$ can be refined to a total ordered $\leqslant$ on G. So, let $s' = \pi(h)$ be the minimal element in *supp h*. Hence, $0 = (hk)(s_0 + s') = h(s_0)\sigma_{s_0}k(s')\tau(s_0,s') = h(s_0)k(s') = bk(s')$, since R is $\sigma$-compatible and $\tau(s_0,s')$ is an invertible element which contradicts the fact that $b \in Z_\ell^*(R)$. Hence $h \in Z_\ell^*(R((G;\sigma;\tau)))$ and we get that $f = g + h$ where $g \in Z_\ell(R((G;\sigma;\tau)))$ and $h \in Z_\ell^*(R((G;\sigma;\tau)))$. Therefore $R((G;\sigma;\tau))$ is left fusible. A right ideal $I$ of a ring $R$ is said to be essential (or large), denoted by $I \leq_e R_R$, if for every right ideal $L$ of $R$, $I \cap L = 0$ implies that $L = 0$. Following [7], the right singular ideal of a ring $R$ is denoted by $Sing(R_R) = \{x \in R | r(x) \leq_e R\}$ where $r(x)$ denotes the right annihilator of $x$. A ring $R$ is called right nonsingular if $Sing(R_R) = 0$. Similarly, we can define $Sing(R((G;\sigma;\tau))_{R((G;\sigma;\tau))})$.

**Corollary 3.3.** Let $(G, \leq)$ a quasitotally ordered group, $\sigma: G \to Aut(R)$ and $\tau: G \times G \to U(R)$, the group of units in $R$. If $R$ is a $\sigma$-compatible and left fusible ring, then $R((G;\sigma;\tau))$ is a right nonsingular ring.

**Proof.** By proposition 2.1, $R((G;\sigma;\tau))$ is a left fusible ring yielding that $R((G;\sigma;\tau))$ is a right nonsingular ring by [7, proposition 2.11].

## 4. SA rings of Mal'cev–Neumann Series.

In this section, we study the right *IN* and right *SA* on rings of Mal'cev–Neumann series.

**Definition 4.1.** A ring $R$ is said to be a right *SA*-ring, if for any two ideals $I, J$ of $R$ there is an ideal $K$ of $R$ such that $r(I) + r(J) = r(K)$, where $r(I)$ denotes the right annihilator of $I$.

**Definition 4.2.** A ring $R$ is called a right *Ikeda–Nakayama* (for short, a right *IN*-ring) if the left annihilator of the intersection of any two right ideals is the sum of the left annihilators; that is, if $\ell(I \cap J) = \ell(I) + \ell(J)$ for all right ideals $I, J$ of $R$; and we say $R$ is an *IN*-ring if $R$ is a left and a right *IN*-ring (see [3],[12]).

**Lemma 4.3.** Let $R$ be a ring, $(G, \leq)$ a quasitotally ordered group and $I$ and $J$ right ideals in $R$. Then $I((G; \sigma; \tau))$ and $J((G; \sigma; \tau))$ are right ideals of $R((G; \sigma; \tau))$ such that $\ell_{R((G;\sigma;\tau))}(I((G; \sigma; \tau)) \cap J((G; \sigma; \tau))) = \ell_{R((G;\sigma;\tau))}((I \cap J)((G; \sigma; \tau))) = \ell_R(I \cap J)((G; \sigma; \tau))$.

**Proof.** It can be easily shown that $I((G; \sigma; \tau)) \cap J((G; \sigma; \tau)) = (I \cap J)((G; \sigma; \tau))$. Hence, using [16, lemma 2.1] the Lemma follows.

**Proposition 4.4.** Let $R$ be a $\sigma$- compatible ring, $(G, \leq)$ a totally ordered group. If $\Lambda = R((G; \sigma; \tau))$ is a right *IN*-ring, then so is $R$.

**Proof.** Let $I$ and $J$ be right ideals of $R$. Then $I((G; \sigma; \tau))$ and $J((G; \sigma; \tau))$ are right ideals of $\Lambda = R((G; \sigma; \tau))$, so by hypothesis, $\ell_\Lambda(I((G; \sigma; \tau)) \cap J((G; \sigma; \tau))) = \ell_\Lambda(I((G; \sigma; \tau))) + \ell_\Lambda(J((G; \sigma; \tau)))$. Hence $\ell_R(I \cap J)((G; \sigma; \tau)) = \ell_\Lambda(I((G; \sigma; \tau))) + \ell_\Lambda(J((G; \sigma; \tau)))$ [16, lemma 2.1]. We need to prove that $\ell_R(I \cap J) = \ell_R(I) + \ell_R(J)$. It is clear that $\ell_R(I) + \ell_R(J) \subseteq \ell_R(I \cap J)$. Now let $a \in \ell_R(I \cap J)$.

Then $f = a\bar{1} \in \ell_\Lambda\big((I \cap J)((G;\sigma;\tau))\big)$. Hence by hypothesis there exists $h = \sum_{y \in G} b_y \bar{y} \in \ell_\Lambda\big(I((G;\sigma;\tau))\big) = \ell_R(I)((G;\sigma;\tau))$ by Lemma 2.3 and $g \in \sum_{z \in G} c_z \bar{z} \in \ell_\Lambda\big(J((G;\sigma;\tau))\big) = \ell_R(J)((G;\sigma;\tau))$ such that $f = a\bar{1} = h + g$. Hence there exist $b$ and $c$ such that $a = b + c$ where $b \in \ell_R(I)$ and $c \in \ell_R(J)$; thus $a \in \ell_R(I) + \ell_R(J)$ and consequently $\ell_R(I \cap J) = \ell_R(I) + \ell_R(J)$.

**Theorem 4.5.** The following statements hold:

(i) If $\Lambda = R((G;\sigma;\tau))$ is a $\sigma$- compatible and right $SA$-ring, then $R$ is a right $SA$-ring;

(ii) If $R$ is a $G$-Armendariz ring, then $R$ is a right $SA$-ring if and only if $\Lambda = R((G;\sigma;\tau))$ is a right $SA$-ring.

**Proof.** (i) Let $I$ and $J$ be right ideals of $R$. Then $I((G;\sigma;\tau))$ and $J((G;\sigma;\tau))$ are right ideals of $R((G;\sigma;\tau))$. So there exist a right ideal $K$ of $R((G;\sigma;\tau))$ such that $r_\Lambda\big(I((G;\sigma;\tau))\big) + r_\Lambda\big(J((G;\sigma;\tau))\big) = r_\Lambda(K)$. Now let $K_0 = \bigcup_{f \in K} C(f)$, then it follows that $K_0$ is a right ideal of $R$. We prove that $r_R(I) + r_R(J) = r_R(K_0)$. Suppose that $a \in r_R(I)$ and $b \in r_R(J)$. Then $f = a\bar{1} \in r_\Lambda\big(I((G;\sigma;\tau))\big) = r_R(I)((G;\sigma;\tau))$ and $g = b\bar{1} \in r_\Lambda\big(J((G;\sigma;\tau))\big) = r_R(J)((G;\sigma;\tau))$, so by hypothesis $f + g \in r_\Lambda(K)$. Then for each $h \in K$, $h(a\bar{1} + b\bar{1}) = 0$. So, $a + b = r(C(h)) \subseteq r_R(K_0)$. Therefore, $r_R(I) + r_R(J) \subseteq r_R(K_0)$. Now let, $c \in r_R(K_0)$. Then, $h = c\bar{1} \in r_\Lambda(K)$. By assumption there exists $f = \sum_{x \in G} a\bar{x} \in r_R(I)((G;\sigma;\tau))$ and $g = \sum_{y \in G} b\bar{y} \in r_R(J)((G;\sigma;\tau))$ such that, $f + g = h$. Hence there exist $a_x \in r_R(I)$ and $b_y \in$ $r_R(J)$ such that $a_x + b_y = c \in r_R(I) + r_R(J)$. So, $r_R(K_0) \subseteq r_R(I) + r_R(J)$. Consequently $r_R(K_0) = r_R(I) + r_R(J)$. Hence, $R$ is a right *SA*-ring.

(ii) The necessity is evident by (i). Now, let $R$ be a $G$-Armendariz and right *SA*-ring and $I$ and $J$ are right two ideals of $R((G; \sigma; \tau))$.

Let, $I_0 = \bigcup_{f \in I} C(f)$ and $J_0 = \bigcup_{f \in J} C(f)$ are right two ideals of $R$, where $C(f)$ denotes the content of $f$. Then, there exists an ideal $K$ in $R$ such that

$r_R(K) = r_R(I_0) + r_R(J_0)$. Now we prove that $r_\Lambda(I) + r_\Lambda(J) = r_\Lambda(K((G; \sigma; \tau)))$. It is sufficient to show that $r_\Lambda(I) + r_\Lambda(J) \subseteq r_\Lambda(K((G; \sigma; \tau)))$. Let $f = \sum_{x \in G} a_x \overline{x} \in r_\Lambda(I)$ and $g = \sum_{y \in G} b_y \overline{y} \in r_\Lambda(J)$, Then for each $h \in I$ and $\rho \in J$, $hf = 0$ and $\rho g = 0$. Since, $R$ is $G$-Armendariz, then $h_k . a_x = 0$ and $\rho_z . b_y = 0$ for all $k \in supp(h)$, $x \in supp(f)$, $z \in supp(\rho)$ and $y \in supp(g)$. Therefore, $a_x \in r(I_0)$ and $b_y \in r(J_0)$. So $a_x + b_y \in r_R(I_0) + r_R(J_0) = r_R(K)$, i.e., there exists $c \in r_R(K)$ such that $a_x + b_y = c$, then for all $m \in G$, $\sum_{m \in G} a_x \overline{m} + \sum_{m \in G} b_y \overline{m} = \sum_{m \in G} c \overline{m}$ which implies that $f + g \in r_\Lambda(K((G; \sigma; \tau)))$. Thus $r_\Lambda(I) + r_\Lambda(J) \subseteq r_\Lambda(K((G; \sigma; \tau)))$. Therefore $r_\Lambda(I) + r_\Lambda(J) = r_\Lambda(K((G; \sigma; \tau)))$.

## 5. Σ-Zip Rings of Mal'cev–Neumann series

In this section, we investigate $\Sigma_U$–zip property in the ring $\Lambda$ of Mal'cev- Neumann series.

**Definition 5.1** [14, Definition 4.1]. Let $\sigma: S \to End(R)$ be a monoid homomorphism and $U$ an ideal of $R$. We say that $U$ is $\Sigma$-compatible if for each $a, b \in R$ and each $s \in S$, $ab \in U \leftrightarrow a\sigma_s(b) \in U$.

**Definition 5.2** [14]. An ideal $U$ of a ring $R$ is called semiprime if for any $a \in R$, $a^n \in U$ implies $a \in U$.

**Lemma 5.3** [14, Lemma 4.2]. Let $\sigma: S \to End(R)$ be a monoid homomorphism and $U$ an ideal of $R$. If $U$ is $\Sigma$-compatible, then for each $a, b \in R$ and each $s \in S$, $ab \in U \leftrightarrow \sigma_s(a)b \in U$.

**Theorem 5.4.** Let $R$ be a ring, $U$ a $\Sigma$-compatible semiprime ideal of $R$ and $G$ a totally ordered group, Then, $R$ is $\Sigma_U$-zip ring if and only if $\Lambda = R((G; \sigma; \tau))$ is a $\Sigma_{U((G; \sigma; \tau))}$-zip ring.

**Proof.** Suppose that $\Lambda = R((G; \sigma; \tau))$ is a $\Sigma_{U((G; \sigma; \tau))}$-zip and $Y \subseteq R$ such that $Y \nsubseteq U$ and $(U:Y) = U$. Let $\bar{Y} = \{y\bar{1} \mid y \in Y\} \subseteq R((G; \sigma; \tau))$. (we need to show that there exists a finite subset $Y_0 \subseteq Y$ such that $(U:Y_0) = U$). It is clear that $U((G; \sigma; \tau)) = \{f \in R((G; \sigma; \tau)) \mid f(y) \in U \text{ for each } y \in supp\ f\}$, is such that $U((G; \sigma; \tau)) \subseteq (U((G; \sigma; \tau)):\bar{Y})$ and it is sufficient to show that $(U((G; \sigma; \tau)):\bar{Y}) \subseteq U((G; \sigma; \tau))$. So, let $f \in (U((G; \sigma; \tau)):\bar{Y})$, therefore $(\bar{y}f)(s) = y\bar{1}f(s) = y\sigma_1 f(s)\tau(1, s) = yf(s)\tau(1, s) \in U$, since $\tau(1, s)$ is an invertible element in $R$ we conclude that $yf(s) \in U$. Hence, $f(s) \in U$ for each $s \in supp\ f$. Consequently $f \in U((G; \sigma; \tau))$ and it follows that $U((G; \sigma; \tau)) = (U((G; \sigma; \tau)):\bar{Y})$

Since $R((G;\sigma;\tau))$ is $\Sigma_{U((G;\sigma;\tau))}$-zip, it follows that there exists a finite subset $\bar{Y}_0 \subseteq \bar{Y}$ such that $(U((G;\sigma;\tau)):\bar{Y}_0) = U((G;\sigma;\tau))$. So, $Y_0 = \{y \in Y \mid y\bar{1} \in \bar{Y}_0\}$ is a finite subset of $Y$ and we get $(U((G;\sigma;\tau)):\bar{Y}_0) \cap R = U((G;\sigma;\tau)) \cap R = U$. Therefore, $R$ is a $\Sigma_U$-zip ring.

Conversely, suppose that $R$ is a $\Sigma_U$-zip ring and $X \subseteq R((G;\sigma;\tau))$ such that $X \nsubseteq U((G;\sigma;\tau))$ and $(U((G;\sigma;\tau)):X) = U((G;\sigma;\tau))$. Let $C_X = \cup_{f \in X} C_f = \cup_{f \in X} \{f(s) \mid s \in supp\, f\} \subseteq R$ be the content of all element of $X$.

We need to show that $(U:C_X) = U$. It is clear that $U \subseteq (U:C_X)$. So, it is sufficient to show that $(U:C_X) \subseteq U$. Let $r \in (U:C_X)$. Then, $ar \in U$ for each $a \in C_X$. Since U is a $\Sigma$-compatible ideal, then $a\sigma_s(r) \in U$ for each $s \in G$. Hence for each $f \in X$ and $s \in G, (fr\bar{1})(s) = f(s)\sigma_s(r)\tau(1,s) \in U$. Therefore, $r\bar{1} \in (U((G;\sigma;\tau)):X) = U((G;\sigma;\tau))$. Hence, $r \in U$ and $(U:C_X) = U$ follows. Since $R$ is a $\Sigma_U$-zip ring, then there exists a finite subset $C_{X_0}$ of $C_X$ such that $(U:\, C_{X_0}) = U$. Let, $X_0 = \{f \in X \mid f(s) \in C_{X_0}$ for some $s \in supp\, f\}$ be a minimal subset of $X$ and it is clear that $X_0$ is finite. We need to show that $(U((G;\sigma;\tau)):X_0) = U((G;\sigma;\tau))$. It is clear that $U((G;\sigma;\tau)) \subseteq (U((G;\sigma;\tau)):X_0)$. So, it is sufficient to show that $(U((G;\sigma;\tau)):X_0) \subseteq U((G;\sigma;\tau))$. So, let $g \in (U((G;\sigma;\tau)):X_0)$, with $v$ be the minimal element of $supp\, g$. Then for each $f \in X_0$, with the minimal element $u$ of $supp\, f$, $(fg)(uv) = f(u)\sigma_u(g(v))\tau(u,v) \in U$. Therefore, $f(u)\sigma_u(g(v))\tau(u,v) = f(u)\sigma_u\left(g(v)\sigma_u^{-1}(\tau(u,v))\right) \in U$. Since, $U$ is a semiprime and Σ-compatible ideal, then $f(u)\left(g(v)\sigma_u^{-1}(\tau(u,v))\right)$ and $\sigma_s\left(g(v)\sigma_u^{-1}(\tau(u,v))\right)f(u) \in U$ for each $s \in G$.

Now, suppose that $w \in G$ is such that for each $u \in supp\ f$ and each $v \in supp\ g$, such that $uv < w$, $\sigma_s\left(g(v)\sigma_u^{-1}(\tau(u,v))\right)f(u) \in U$ for each $s \in G$. Using transfinite induction, we need to show that $f(u_i)\sigma_{u_i}g(v_i)\tau(u_i,v_i) \in U$ for each $uv = w$. Since, $X_w(f,g) = \{(u,v) \in G \times G | uv = w, u \in suppf$ and $v \in supp\ g\}$ is a finite set and $G$ totally ordered then, let $\{u_i, v_i, i = 1,2,3,\ldots,n\}$ be such that $u_1 < u_2 \ldots \ldots < u_n$ and $v_n < v_{n-1} \ldots \ldots < v_1$. Hence, $(fg)(w) = \sum_{\substack{(u,v)\in X_w(f,g) \\ i=1}}^{n} f(u_i)\sigma_{u_i}(g(v_i))\tau(u_i,v_i) =$

$f(u_1)\sigma_{u_1}(g(v_1))\tau(u_1,v_1) + \cdots + f(u_n)\sigma_{u_n}(g(v_n))\tau(u_n,v_n) = a_1 \in U$ (1)

Since, for each $u_i, i \geq 2$, $u_1 v_i < u_i v_i = w$, then using induction hypothesis, we have $\sigma_{u_i}\left(g(v_i)\sigma_{u_i}^{-1}(\tau(u_i,v_i))\right)f(u_1) \in U$. Multiplying (1) on the right by $f(u_1)$, then we have $fg(w)f(u_1) = f(u_1)\sigma_{u_1}(g(v_1))\ \tau(u_1,v_1)f(u_1) +$ $f(u_2)\sigma_{u_2}\left(g(v_2)\sigma_{u_2}^{-1}(\tau(u_2,v_2))\right)f(u_1) + \cdots +$ $f(u_n)\sigma_{u_n}\left(g(v_n)\sigma_{u_n}^{-1}(\tau(u_n,v_n))\right)f(u_1) = a_1 f(u_1) \in U$. Thus, we obtain $f(u_1)\sigma_{u_1}\left(g(v_1)\sigma_{u_1}^{-1}(\tau(u_1,v_1))\right)f(u_1) \in U$. Since $U$ is semiprime it follows that $f(u_1)\sigma_{u_1}(g(v_1))\ \tau(u_1,v_1) \in U$. Now, subtract $f(u_1)\sigma_{u_1}(g(v_1))\ \tau(u_1,v_1)$ and multiply by $f(u_2)$ from the right of both sides of (1), it follows that $(fg)(w)f(u_2) - f(u_1)\sigma_{u_1}(g(v_1))\ \tau(u_1,v_1)f(u_2) =$ $f(u_2)\sigma_{u_2}(g(v_2))\ \tau(u_2,v_2)f(u_2) + \cdots + f(u_n)\sigma_{u_n}(g(v_n))\ \tau(u_n,v_n)f(u_2)=$ $a_1 f(u_2) - f(u_1)\sigma_{u_1}(g(v_1))\,\tau(u_1,v_1)f(u_2) \in U$. Using the same argument as above we obtain $f(u_2)\sigma_{u_2}(g(v_2))\,\tau(u_2,v_2) \in U$. Continuing this process, we can show that $f(u_i)\sigma_{u_i}(g(v_i))\,\tau(u_i,v_i) \in U$ for each $i = 1,2,3,\ldots,n$ such that $u_i v_i = w$. Thus, $f(u)\sigma_u(g(v))\tau(u,v) \in U$ for each

$u \in supp f$ and $v \in supp g$. Hence $g \in U((G;\sigma;\tau))$ and it follows that $(U((G;\sigma;\tau)):X_0) \subseteq U((G;\sigma;\tau))$. Consequently, we deduce that $R((G;\sigma;\tau))$ is a $\Sigma_{U((G;\sigma;\tau))}$-zip ring.

In the following we give some examples of $\Sigma_U$- zip ring

**Example 5.5.** Let $R = Z_4$ be the ring of integer modulo 4 and $U = <2>$ the ideal generated by 2, then $R$ is a $\Sigma_U$- zip ring since $(U:X) = U$ for each subset $X \nsubseteq U$. If $U = \{0\}$, then $R$ is $\Sigma_0$- zip as well as zip.

**Example 5.6.** Let $R = T(Z_4, Z_4) \cong \left\{ \begin{pmatrix} a & b \\ 0 & a \end{pmatrix}, a,b \in Z_4 \right\}$ the trivial extension of $Z_4$. We can write the proper ideal of $R$ as $U = \left\{ \begin{pmatrix} 0 & m \\ 0 & 0 \end{pmatrix} : m \in Z_4 \right\}$ let $X$ be any subset of $R$ with $X \nsubseteq U$, then $R$ is a $\Sigma_U$- zip ring since $(U:X) = U$ and for any subset $X_0$ of $X$, we have $(U:X_0) = U$ by a routine computation. So $R$ is $\Sigma_U$- zip.

## References

[1] Beachy J. A. and Blair W. D. (1975). Rings whose faithful left ideals are cofaithful. *Pacific J. Math.* 58, 1–13.

[2] Birkenmeier G. F., Ghirati M. and Taherifar A. (2015).When is a Sum of Annihilator Ideals an Annihilator Ideal, Comm. Algebra 43, 2690-2702.

[3] Camillo V., Nicholson W. K. and Yousif M. F. (2000). Ikeda- Nakayama rings. J. Algebra 266, 1001-1010.

[4] Cedo F. (1991).Zip rings and Mal'cev domains. Comm. Algebra 19, 1983–1991.

[5] Cortes W. (2008). Skew polynomial extensions over zip rings. Int. J. Math. Sci. 10, 1–8.

[6] Faith C. (1989). Rings with zero intersection property on annihilators: zip rings. Publ. Math. 33, 329–332.

[7] Ghashghaei E. and Mcgovern W. W. (2017). Fusible rings. Comm. Algebra 45, 1151–1165.

[8] Hong C. Y., Kim N. K., Kwak T. K. and Lee Y. (2005).Extension of zip rings. J. Pure Appl. Algebra 195, 231–242.

[9] Lorenz E. M. (1983). Division algebras generated by finitely generated nilpotent groups. J. Algebra 85,368–381.

[10] Makar L.-Limanov(1983). The skew field of fractions of the Weyl algebra contains a free noncommutative subalgebra, Comm. Algebra 11, 2003–2006.

[11] Musson I. M. and Stafford K. (1993). Mal'cev-Neumann group rings. Comm. Algebra 21, 2065–2075.

[12] Nicholson W. K. and Yousif M. F. (2003). Quasi-Frobenius Rings, Cambridge University Press.

[13] Ouyang L.(2009). Ore extensions of weak zip rings, Glasg. Math. J. 51, 525–537.

[14] Ouyang L., Zhou Q. and Jinfang Wu. (2017). Extensions of Σ-Zip rings. Int. Elec. J. Algebra 21,1-22.

[15] Salem R. M. (2013). Generalized power series over zip and weak zip rings. Southeast Asian Bull. Math. 37,259–268.

[16] Salem R. M. (2012).Mal'cev-Neumann series over zip and Weak zip rings. Asian- European J. Math. 4,1250058-1250066.

[17] Sonin C. (1998). Krull dimension of Mal'cev-Neumann rings, Comm. Algebra 26,2915–2931.

[18] Zelmanowitz J. M. (1976).The finite intersection property on annihilator right ideals. Proc. Amer. Math. Soc. 57,213–216.

[19] Zhao R.and Liu Z. (2008). On some properties of Mal'cev-Neumann modules. Bull. Korean Math. Soc. 45,445–456.